\documentclass[english,11pt, a4paper]{article}

\title{Ahlfors $1$-Regularity and Conformal Criteria
for Carleson Measures}

\usepackage{iftex}
\usepackage[T1]{fontenc}

\usepackage[utf8]{inputenc}

\usepackage{mathtools}
\usepackage{lmodern}
\usepackage{amsthm,amsmath,amsfonts,amssymb}
\usepackage{graphicx}
\usepackage{xcolor}

\definecolor{revisionred}{RGB}{150,35,35}

\usepackage{tikz}
\usetikzlibrary{arrows.meta,calc}
\usepackage{enumerate,enumitem}
\usepackage{authblk}
\usepackage{cite}

\let\OLDthebibliography\thebibliography
\renewcommand\thebibliography[1]{
  \OLDthebibliography{#1}
  \setlength{\parskip}{0pt}
  \setlength{\itemsep}{0pt plus 0.3ex}
}

\usepackage[activate={true,nocompatibility},final,tracking=true,kerning=true,spacing=false,factor=1100,stretch=15,shrink=15]{microtype}
\ifXeTeX
  \microtypesetup{expansion=false,tracking=false,kerning=false}
\fi

\allowdisplaybreaks

\numberwithin{equation}{section}

\theoremstyle{plain}
\newtheorem{theorem}{Theorem}

\newtheorem{proposition}[theorem]{Proposition}

\theoremstyle{remark}

\usepackage{babel}
\usepackage[hypertexnames=false]{hyperref}
\hypersetup{
    colorlinks=true,
    allcolors=black, 
    linktoc=all 
}

\author{%
Yihao Hou\thanks{Center for applied Mathematics, Tianjin University,
Weijin Road 92, Tianjin, 300072, P.R. China, email: \url{yhhou@tju.edu.cn}}
\, and Huaying Wei\thanks{Center for applied Mathematics and KL-AAGDM,
Tianjin University, Weijin Road 92, Tianjin, 300072, P.R. China,
email: \url{hywei@tju.edu.cn}}}
\date{}

\newcommand{\HH}{\mathcal H^1}
\newcommand{\CM}{\mathrm{CM}}

\newcommand{\length}{\operatorname{length}}

\newcommand{\dist}{\operatorname{dist}}

\begin{document}

\maketitle

\begin{abstract}
We prove that a bounded Jordan domain $\Omega$ has
Ahlfors $1$-regular boundary if and only if the following
criterion holds for every positive Borel measure $\mu$
on $\Omega$: $\mu$ is Carleson exactly when
$\int_\Omega |\phi'(w)|\,d\mu(w)<\infty$ for every conformal
bijection $\phi$ from $\Omega$ onto a bounded domain with
rectifiable boundary.
The same characterization holds using only arclength
measures on simple locally rectifiable curves, with the
integral condition replaced by finite length of every
such conformal image.
These results extend Bishop's theorem on the unit disk.
We also study analogous criteria for vanishing Carleson
measures.


\end{abstract}


\section{Introduction and main results}

Let $\Omega$ be a bounded Jordan domain with boundary
$\Gamma=\partial\Omega$. A positive Borel measure $\mu$ on
$\Omega$ is called a \emph{Carleson measure}, written
$\mu\in\CM(\Omega)$, if
\[
 \|\mu\|_{\CM(\Omega)}
 :=\sup_{x\in\Gamma,\,r>0}
 \frac{\mu(B(x,r)\cap\Omega)}{r}<\infty, 
\]
where $B(x, r)$ is the disk with center $x$ and radius $r$. 
Moreover, 
a measure $\mu\in\CM(\Omega)$ is called a
\emph{vanishing Carleson measure}, written
$\mu\in\CM_0(\Omega)$, if
\begin{equation}\label{eq:vanishing-definition}
 \lim_{r\downarrow0}
 \sup_{x\in\Gamma}
 \frac{\mu(B(x,r)\cap\Omega)}{r}=0.
\end{equation}

Bishop~\cite{Bishop} proved that a positive Borel measure
$\mu$ on the unit disk $\mathbb D$ is a Carleson measure
if and only if
\[
 \int_{\mathbb D}|G'(z)|\,d\mu(z)<\infty
\]
for every conformal map $G$ from $\mathbb D$ onto a bounded
domain with rectifiable boundary. In particular, when
$\mu$ is arclength measure on a simple locally rectifiable
curve $\gamma\subset\mathbb D$, this criterion says that
arclength measure on $\gamma$ is Carleson if and only if $G(\gamma)$
has finite length for every such map $G$.

In this note, we characterize Ahlfors $1$-regularity of
Jordan curves by extending Bishop's criterion to bounded
Jordan domains, thereby resolving a conjecture proposed
by Michel Zinsmeister (personal communication).
Recall that a Jordan curve $\Gamma$ is
\emph{Ahlfors $1$-regular} if there exists a constant $C>0$
such that
\[
 \HH(\Gamma\cap B(x,r))\le Cr,
 \qquad x\in\Gamma,\quad
 0<r<\operatorname{diam}\Gamma.
\]

We call a conformal bijection $\phi:\Omega\to\Omega'$
\emph{admissible} if $\Omega'$ is a bounded domain with
rectifiable boundary. For a simple locally rectifiable
curve $\gamma\subset\Omega$, write
$\sigma_\gamma=\HH\!\restriction_\gamma$ for its arclength
measure.

\begin{theorem}\label{thm:measures}
Let $\Omega$ be a bounded Jordan domain.  The following statements are equivalent:
\begin{enumerate}
\item[\textup{($C_1$)}] $\Gamma :=\partial\Omega$ is Ahlfors $1$-regular.

\item[\textup{($C_2$)}] For every positive Borel measure
$\mu$ on $\Omega$,
\begin{equation}\label{eq:measure-test}
 \mu\in\CM(\Omega)
 \quad\Longleftrightarrow\quad
 \int_\Omega|\phi'(w)|\,d\mu(w)<\infty
 \quad\text{for every admissible }\phi.
\end{equation}

\item[\textup{($C_3$)}] For every simple locally rectifiable
curve $\gamma\subset\Omega$,
\[
 \sigma_\gamma\in\CM(\Omega)
 \quad\Longleftrightarrow\quad
 \length(\phi(\gamma))<\infty
 \quad\text{for every admissible }\phi.
\]
\end{enumerate}
\end{theorem}

The implication $\textup{$C_2$}\Rightarrow\textup{$C_3$}$
follows by taking $\mu=\sigma_\gamma$. The theorem shows
that restricting the testing measures to arclength measures
does not weaken the characterization.

Fix a Riemann map $F:\mathbb D\to\Omega$, and define
\[
 m_a(z)=\frac{a-z}{1-\overline a z},
 \qquad
 \phi_a=m_a\circ F^{-1},
 \qquad a\in\mathbb D.
\]
Each $\phi_a$ is an admissible conformal map onto $\mathbb D$.
Whereas Theorem~\ref{thm:measures} tests finiteness against
every admissible map, the following vanishing counterpart
tests decay along the canonical family
$(\phi_a)_{a\in\mathbb D}$ as $|a|\to1$.

\begin{theorem}
\label{thm:vanishing-positive-measures}
Let $\Omega$ be a bounded Jordan domain, and let
$(\phi_a)_{a\in\mathbb D}$ be defined as above.
The following statements are equivalent:
\begin{enumerate}
\item[\textup{($V_1$)}] $\Gamma := \partial\Omega$ is Ahlfors
$1$-regular.

\item[\textup{($V_2$)}] For every positive Borel measure
$\mu$ on $\Omega$,
\begin{equation}\label{eq:vanishing-test}
 \mu\in\CM_0(\Omega)
 \quad\Longleftrightarrow\quad
 \lim_{|a|\to1}
 \int_\Omega|\phi_a'(w)|\,d\mu(w)=0.
\end{equation}
Here, the limit is uniform in the argument of $a$.
\end{enumerate}
\end{theorem}

The restriction to arclength measures behaves differently
in the vanishing setting. Indeed, the natural curve analogue
of \textup{$V_2$} holds on every bounded Jordan domain and
therefore cannot characterize Ahlfors $1$-regularity as
\textup{$C_3$} does. The following proposition makes this
distinction precise.

\begin{proposition}\label{prop:non}
Let $\Omega$ be a bounded Jordan domain, and let
$(\phi_a)_{a\in\mathbb D}$ be defined as above.
For every simple locally rectifiable curve
$\gamma\subset\Omega$,
\[
 \sigma_\gamma\in\CM_0(\Omega)
 \quad\Longleftrightarrow\quad
 \lim_{|a|\to1}
 \int_\gamma|\phi_a'(w)|\,ds(w)=0,
\]
where $ds$ means the element of arclength along the curve $\gamma$, and  the limit is uniform in the argument of $a$.
No Ahlfors $1$-regularity assumption on $\Gamma :=\partial\Omega$
is required.
\end{proposition}

The remainder of the note is organized as follows.
Section~2 is devoted to Theorem~\ref{thm:measures}, while
Section~3 establishes Theorem~\ref{thm:vanishing-positive-measures}
and Proposition~\ref{prop:non}.

\section{Carleson measures and conformal images: Theorem~\ref{thm:measures}}
\label{sec:carleson-images}

\subsection{Carleson measure criteria}

Let $F:\mathbb D\to\Omega$ be a Riemann map. For a positive Borel measure $\mu$ on $\Omega$,
define its weighted pullback by
\begin{equation}\label{eq:pullback}
 F^*\mu(A)
 :=
 \int_{F(A)}|(F^{-1})'(w)|\,d\mu(w),
 \qquad A\subset\mathbb D\text{ Borel}.
\end{equation}
In particular, the weighted pullback transforms the arclength measure $\sigma_\gamma$ on a curve $\gamma$ into the arclength measure $\sigma_{F^{-1}(\gamma)}$ on its conformal preimage $F^{-1}(\gamma)$. 
If
$d\mu(w)=\lambda(w)\,dA(w)$, the change-of-variables formula
gives
\[
 d(F^*\mu)(z)=|F'(z)|^{-1}d\mu(F(z)) = \lambda(F(z))|F'(z)|\,dA(z),
\]
where $dA$ denotes the planar area measure. 
Conversely, for a positive Borel measure $\nu$ on $\mathbb D$,
define its weighted pushforward by
\begin{equation}\label{eq:pushforward}
 (F^{-1})^*\nu(E)
 :=
 \int_{F^{-1}(E)}|F'(z)|\,d\nu(z),
 \qquad E\subset\Omega\text{ Borel}.
\end{equation}
These two operations are inverse to one another.  For an admissible map $\phi:\Omega\to\Omega'$, put
$G=\phi\circ F$. The maps and measures are related by the following commutative diagram: 
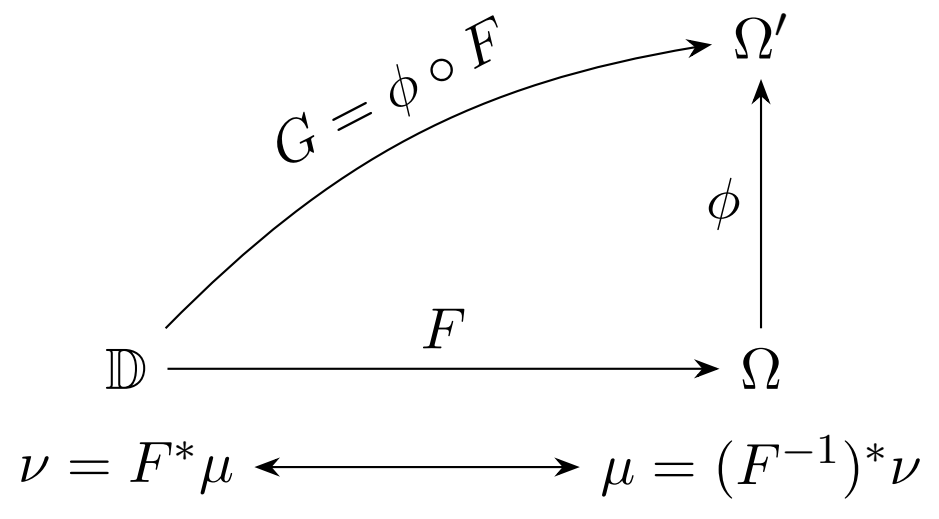
\[
\begin{tikzpicture}[
    >=Stealth,
    baseline=(current bounding box.center)
]
 \node (d) at (0,0) {$\mathbb D$};
 \node (o) at (4.2,0) {$\Omega$};
 \node (p) at (4.2,2.2) {$\Omega'$};

 \draw[->] (d) -- node[above] {$F$} (o);
 \draw[->] (o) -- node[left] {$\phi$} (p);
 \draw[->] (d) to[bend left=18]
     node[above,sloped] {$G=\phi\circ F$} (p);

 \node (nu) at (0,-0.65) {$\nu=F^*\mu$};
 \node (mu) at (4.2,-0.65) {$\mu=(F^{-1})^*\nu$};
 \draw[<->] (nu.east) -- (mu.west);
\end{tikzpicture}
\]

We use the following three measure criteria.
\begin{enumerate}[leftmargin=*,itemsep=0.5em]
\item[] \textbf{Bishop's criterion.}
As we have recalled above, for every positive Borel measure $\nu$ on $\mathbb D$,
\begin{equation}\label{eq:bishop}
 \nu\in\CM(\mathbb D)
 \quad\Longleftrightarrow\quad
 \int_{\mathbb D}|G'(z)|\,d\nu(z)<\infty
 \quad\text{for every admissible }G.
\end{equation}

\item[] \textbf{Transfer on Ahlfors $1$-regular domains.}
If $\Gamma=\partial\Omega$ is Ahlfors $1$-regular, both
weighted transfers are bounded on the corresponding
Carleson spaces. In particular,
\begin{equation}\label{eq:transfer}
 \mu\in\CM(\Omega)
 \quad\Longleftrightarrow\quad
 F^*\mu\in\CM(\mathbb D).
\end{equation}
These transfer results go back to
Zinsmeister~\cite{Zinsmeister}; see also
Wei-Zinsmeister~\cite{WeiZinsmeister}.

\item[] \textbf{The pushforward characterization.}
Wei-Zinsmeister~\cite{WeiZinsmeister} proved that
$\Gamma$ is Ahlfors $1$-regular if and only if
\begin{equation}\label{eq:push-characterization}
 (F^{-1})^*\nu\in\CM(\Omega)
 \quad\text{for every }\nu\in\CM(\mathbb D).
\end{equation}
\end{enumerate}

\begin{proof}[Proof of Theorem \ref{thm:measures}]
We prove $\textup{$C_1$}\Leftrightarrow
\textup{$C_2$}$ and
$\textup{$C_2$}\Rightarrow\textup{$C_3$}\Rightarrow\textup{$C_1$}$. 
Given a positive Borel measure $\mu$ on $\Omega$, put
$\nu=F^*\mu$. For every admissible $\phi$, 
applying \eqref{eq:pullback} gives 
\begin{equation}\label{eq:identity}
\begin{aligned}
 \int_{\mathbb D}|G'(z)|\,d\nu(z)
 &=
 \int_\Omega
 |G'(F^{-1}(w))|\,|(F^{-1})'(w)|\,d\mu(w)\\
 &=\int_\Omega|\phi'(w)|\,d\mu(w),
 \qquad G=\phi\circ F.
\end{aligned}
\end{equation}
As $\phi$ ranges over the admissible maps on $\Omega$,
$G$ ranges over the admissible maps on $\mathbb D$.

Under \textup{$C_1$}, combining \eqref{eq:transfer},
\eqref{eq:bishop}, and \eqref{eq:identity} therefore proves
\textup{$C_2$}. Conversely, assume \textup{$C_2$}. Given any positive
$\nu\in\CM(\mathbb D)$, set $\mu=(F^{-1})^*\nu$.
Then $F^*\mu=\nu$, so Bishop's criterion and
\eqref{eq:identity} yield
\[
 \int_\Omega|\phi'(w)|\,d\mu(w)<\infty
 \qquad\text{for every admissible }\phi.
\]
Condition \textup{$C_2$} implies $\mu\in\CM(\Omega)$.
Since $\nu$ was arbitrary, \eqref{eq:push-characterization}
gives \textup{$C_1$}. Taking $\mu=\sigma_\gamma$ in \textup{$C_2$}
and using
\begin{equation}\label{eq:curve-length}
 \int_\Omega|\phi'|\,d\sigma_\gamma
 =
 \int_\gamma|\phi'(w)|\,ds(w)
 =
 \length(\phi(\gamma)), 
\end{equation}
we get $\textup{$C_2$}\Rightarrow\textup{$C_3$}$. These directions are easy consequences of known facts, while \textup{$C_3$} $\Rightarrow$
\textup{$C_1$} requires more work. 

We now prove \textup{$C_3$} $\Rightarrow$
\textup{$C_1$} by contradiction. Suppose that
$\Gamma=\partial\Omega$ is not Ahlfors $1$-regular. The central idea of this proof is to turn the failure of Ahlfors $1$-regularity into one curve that contradicts 
\(\textup{$C_3$}\).

\noindent\textbf{Step 1: Localization of bad circular arcs.}
Fix
a Riemann map $F:\mathbb D\to\Omega$. 
We first observe that a uniform bound
\begin{equation}\label{eq:uniform}
    \sup_{1/2<r<1}
 \|\sigma_{F(r\mathbb T)}\|_{\CM(\Omega)}
 \le M<\infty
\end{equation} 
would imply that $\Gamma$ is Ahlfors $1$-regular.

We now give an explicit explanation for this observation.
Write 
\[
L(r):=\operatorname{length}(F(r\mathbb T))
=r\int_0^{2\pi}|F'(re^{it})|\,dt, \quad 0<r<1.
\]
Since \(|F'|\) is subharmonic, its circular means are nondecreasing. Thus, for \(0<r<s<1\),
 $L(r)\le r/sL(s)\le L(s)$ so that \(L(r)\) is nondecreasing. Fix $x_0\in\Gamma$ and $R\geq\operatorname{diam}\Omega$ we   have 
$$L(r)
 =\sigma_{F(r\mathbb T)}(B(x_0,R))
 \le MR$$
 by \eqref{eq:uniform}. 
On the other hand,  Carath\'eodory's theorem gives $F(re^{it})\to F(e^{it})$ uniformly with respect to $t\in [0, 2\pi)$ as $r\uparrow1$.  
Thus, we get 
\[
\operatorname{length}(\Gamma)
\le \liminf_{r\uparrow1}L(r)
=\lim_{r\uparrow1}L(r)
\le MR<\infty,
\]
where the first inequality follows from the definition of curve length using polygonal sums, and the equality follows from the  monotonicity of \(L(r)\). 
This implies that $\Gamma$ is rectifiable, and in particular, $F'\in H^1(\mathbb D)$. Fix $x\in\Gamma$ and $h>0$.  The set
\[
U:=\{t\in(0,2\pi):F(e^{it})\in B(x,h)\}
\]
is open. Write \(U=\bigcup_j I_j\) as the union of its pairwise disjoint open interval components.
We first  consider finitely many intervals $(I_j)_{j=1}^{N}$.  For $r$ sufficiently
close to $1$, the images of these intervals under $F(re^{it})$ lie in
$B(x,2h)$, so their combined length is at most $2Mh$. 
Passing to the limit and applying $F'\in H^1(\mathbb D)$  give
\[
 \sum_{j=1}^{N}\ell\left(\{F(e^{it}):\;t\in I_j\}\right)=
 \lim_{r\uparrow1}
 \sum_{j=1}^{N}\ell\left(\{F(re^{it}):\;t\in I_j\}\right)
 \le2Mh.
\]
Taking $N\to\infty$ we obtain
\[
 \mathcal H^1(\Gamma\cap B(x,h))
 = \lim_{N\to\infty}\sum_{j=1}^{N}\ell\left(\{F(e^{it}):\;t\in I_j\}\right)
 \le2Mh.
\]
This proves the observation.

\bigskip

For $\xi\in\mathbb T$ and $\varepsilon>0$, write
\[
 A(r,\xi,\varepsilon)
 :=
 \{r\xi e^{it}:|t|\le\varepsilon\}.
\]
There consequently exists $\xi\in\mathbb T$ such that
\begin{equation}\label{eq:121}
    \sup_{1/2<r<1}
 \|\sigma_{F(A(r,\xi,\varepsilon))}\|_{\CM(\Omega)}
 =\infty
\end{equation} 
for every sufficiently small $\varepsilon>0$.
Otherwise, each point of $\mathbb T$ would have an
angular neighborhood on which this supremum is finite.
A finite subcover would then give a uniform Carleson
bound for the whole image curves $F(r\mathbb T)$, $1/2<r<1$, contradicting the
preceding observation \eqref{eq:uniform}.

After precomposing $F$ with a rotation, assume that
$\xi=1$. Notice also that, for every $r_0<1$, the image
arcs with $1/2<r\le r_0$ have uniformly bounded lengths
and remain in a fixed compact subset of $\Omega$.
Their Carleson norms are therefore uniformly bounded.
Thus, the bad radii in \eqref{eq:121} can be chosen arbitrarily
close to $1$.

\noindent\textbf{Step 2: Construction of one simple curve.}
Choose successively positive numbers
$\varepsilon_n,\delta_n$ such that
\[
 \varepsilon_1<\frac1{10},
 \qquad
 \delta_n<\frac{\varepsilon_n}{10},
 \qquad
 \varepsilon_{n+1}<\frac{\delta_n}{10}.
\]
Set 
\[
 r_n:=1-\delta_n,
 \qquad
 A_n:=\{r_ne^{it}:|t|\le\varepsilon_n\}. 
\]
We have
\[
 \|\sigma_{F(A_n)}\|_{\CM(\Omega)}>n.
\]
The unboundedness in \eqref{eq:121} makes these choices achievable.

Put
\[
 \rho_n:=1-\frac{\delta_n}{2},
 \qquad s_n:=(-1)^{n+1}.
\]
The endpoints of $A_n$ have angles $s_n\varepsilon_n$. Traverse the arcs $A_n$ in alternating directions and 
connect the endpoint to the starting endpoint of
$A_{n+1}$ by the following:
\begin{enumerate}
\item the radial segment from
$r_ne^{is_n\varepsilon_n}$ to
$\rho_ne^{is_n\varepsilon_n}$;
\item the circular arc at radius $\rho_n$ from angle
$s_n\varepsilon_n$ to $s_n\varepsilon_{n+1}$;
\item the radial segment from
$\rho_ne^{is_n\varepsilon_{n+1}}$ to
$r_{n+1}e^{is_n\varepsilon_{n+1}}$.
\end{enumerate}
Denote the resulting curve by $\eta$; see Figure~\ref{fig:eta-circular-block}.

The scale separation gives
\[
 r_n<\rho_n<r_{n+1}.
\]
Each connector therefore lies in the annulus between
successive main arcs $A_n$. The pieces meet only at consecutive
endpoints, so $\eta$ is simple. It approaches $1$, which
is not included in $\eta$, and every compact subset of
$\mathbb D$ meets only finitely many pieces. In particular,
$\eta$ is locally rectifiable.

\noindent\textbf{Step 3: The curve $\eta$ is Carleson.}
Fix \(\zeta\in\mathbb T\) and \(h>0\). Every point \(z\in B(\zeta,h)\cap\mathbb D\) satisfies
\[
1-|z|\le|\zeta-z|<h.
\]
Along the radial pieces of \(\eta\), the radius increases monotonically, and no radial interval is traversed more than once. Consequently, the total length of the radial pieces contained in \(B(\zeta,h)\) is at most \(h\).

For each \(n\), group the main arc \(A_n\) with the circular part \(C_n\) of its outgoing connector. These two arcs have combined length at most \(3\varepsilon_n\), and their distances from \(\mathbb T\) are \(\delta_n\) and \(\delta_n/2\), respectively. The right-hand panel of Figure~\ref{fig:eta-circular-block}  illustrates this grouping for \(n=4\).  

\begin{figure}[!htb]
\centering
\begingroup
\definecolor{etaBlockBlue}{HTML}{175B94}
\definecolor{etaBlockOrange}{HTML}{C16A2F}
\definecolor{etaBlockGray}{HTML}{929BA3}
\resizebox{0.82\linewidth}{!}{%
\begin{tikzpicture}[
  x=.01cm,y=-.01cm,font=\small,
  line cap=round,line join=round,
  >={Stealth[length=1.7mm,width=1.3mm]},
  etaMain/.style={etaBlockBlue,line width=1.1pt},
  etaConnector/.style={etaBlockOrange,line width=1.1pt},
  etaRadial/.style={etaBlockGray,densely dashed,line width=.6pt},
  etaGap/.style={<->,black,line width=.5pt},
  declare function={
    etazx(\r,\t)=1260+5000*(\r*cos(\t)-1);
    etazy(\r,\t)=440-5000*\r*sin(\t);
  }
]
\path[use as bounding box] (50,100) rectangle (1460,1020);

\draw[line width=.7pt] (300,400) circle[radius=2.2cm];
\foreach \k in {0,...,8}{
  \pgfmathsetmacro{\EtaR}{1-.42*pow(.4,\k)}
  \pgfmathsetmacro{\EtaRho}{1-.21*pow(.4,\k)}
  \pgfmathsetmacro{\EtaNextR}{1-.42*pow(.4,\k+1)}
  \pgfmathsetmacro{\EtaE}{.95*pow(.4,\k)*180/pi}
  \pgfmathsetmacro{\EtaNextE}{.4*\EtaE}
  \pgfmathsetmacro{\EtaS}{pow(-1,\k)}
  \draw[etaBlockBlue,line width=.7pt]
    plot[variable=\t,samples=65,
      domain={-\EtaS*\EtaE}:{\EtaS*\EtaE}]
    ({300+220*\EtaR*cos(\t)},{400-220*\EtaR*sin(\t)});
  \ifnum\k<8
    \draw[etaBlockOrange,line width=.65pt]
      ({300+220*\EtaR*cos(\EtaS*\EtaE)},
       {400-220*\EtaR*sin(\EtaS*\EtaE)}) --
      ({300+220*\EtaRho*cos(\EtaS*\EtaE)},
       {400-220*\EtaRho*sin(\EtaS*\EtaE)});
    \draw[etaBlockOrange,line width=.65pt]
      plot[variable=\t,samples=45,
        domain={\EtaS*\EtaE}:{\EtaS*\EtaNextE}]
      ({300+220*\EtaRho*cos(\t)},{400-220*\EtaRho*sin(\t)});
    \draw[etaBlockOrange,line width=.65pt]
      ({300+220*\EtaRho*cos(\EtaS*\EtaNextE)},
       {400-220*\EtaRho*sin(\EtaS*\EtaNextE)}) --
      ({300+220*\EtaNextR*cos(\EtaS*\EtaNextE)},
       {400-220*\EtaNextR*sin(\EtaS*\EtaNextE)});
  \fi
}
\fill (300,400) circle[radius=.026cm];
\node[below] at (300,408) {$0$};
\node[font=\Large] at (180,402) {$\mathbb D$};
\node at (162,194) {$\mathbb T$};
\fill[etaBlockBlue] (374.2,503.8) circle[radius=.04cm];
\draw[fill=white,line width=.6pt] (520,400) circle[radius=.055cm];
\node[right] at (535,400) {$1$};
\draw[etaBlockGray,densely dashed] (500,374) rectangle (538,426);
\node[text=etaBlockBlue] at (600,325) {$A_4$};
\draw[etaBlockBlue,thin] (582,336)--(514,395);
\draw[->,line width=.8pt] (645,440)--(815,440);
\node[above] at (730,421) {enlarge};

\begin{scope}
\clip (850,100) rectangle (1460,820);
\draw[line width=.75pt]
  plot[variable=\t,domain=-5:5,samples=90]
  ({etazx(1,\t)},{etazy(1,\t)});
\foreach \k in {3,...,8}{
  \pgfmathsetmacro{\EtaR}{1-.42*pow(.4,\k)}
  \pgfmathsetmacro{\EtaRho}{1-.21*pow(.4,\k)}
  \pgfmathsetmacro{\EtaNextR}{1-.42*pow(.4,\k+1)}
  \pgfmathsetmacro{\EtaE}{.95*pow(.4,\k)*180/pi}
  \pgfmathsetmacro{\EtaNextE}{.4*\EtaE}
  \pgfmathsetmacro{\EtaS}{pow(-1,\k)}
  \ifnum\k=3
    \tikzset{etaPiece/.style={etaMain},etaJoin/.style={etaRadial},
      etaMiddle/.style={etaConnector}}
  \else
    \tikzset{etaPiece/.style={gray!35,line width=.55pt},
      etaJoin/.style={gray!35,line width=.55pt},
      etaMiddle/.style={gray!35,line width=.55pt}}
  \fi
  \draw[etaPiece]
    plot[variable=\t,samples=65,
      domain={-\EtaS*\EtaE}:{\EtaS*\EtaE}]
    ({etazx(\EtaR,\t)},{etazy(\EtaR,\t)});
  \ifnum\k<8
    \draw[etaJoin]
      ({etazx(\EtaR,\EtaS*\EtaE)},
       {etazy(\EtaR,\EtaS*\EtaE)}) --
      ({etazx(\EtaRho,\EtaS*\EtaE)},
       {etazy(\EtaRho,\EtaS*\EtaE)});
    \draw[etaMiddle]
      plot[variable=\t,samples=45,
        domain={\EtaS*\EtaE}:{\EtaS*\EtaNextE}]
      ({etazx(\EtaRho,\t)},{etazy(\EtaRho,\t)});
    \draw[etaJoin]
      ({etazx(\EtaRho,\EtaS*\EtaNextE)},
       {etazy(\EtaRho,\EtaS*\EtaNextE)}) --
      ({etazx(\EtaNextR,\EtaS*\EtaNextE)},
       {etazy(\EtaNextR,\EtaS*\EtaNextE)});
  \fi
}

\pgfmathsetmacro{\EtaR}{1-.42*pow(.4,3)}
\pgfmathsetmacro{\EtaRho}{1-.21*pow(.4,3)}
\pgfmathsetmacro{\EtaE}{.95*pow(.4,3)*180/pi}
\pgfmathsetmacro{\EtaAngleA}{\EtaE/2}
\pgfmathsetmacro{\EtaAngleC}{-.7*\EtaE}
\draw[etaGap]
  ({etazx(\EtaR,\EtaAngleA)},{etazy(\EtaR,\EtaAngleA)}) --
  ({etazx(1,\EtaAngleA)},{etazy(1,\EtaAngleA)});
\draw[etaGap]
  ({etazx(\EtaRho,\EtaAngleC)},{etazy(\EtaRho,\EtaAngleC)}) --
  ({etazx(1,\EtaAngleC)},{etazy(1,\EtaAngleC)});
\draw[etaMain,->]
  plot[variable=\t,domain=.86:-.52,samples=20]
  ({etazx(\EtaR,\t)},{etazy(\EtaR,\t)});
\draw[etaConnector,->]
  plot[variable=\t,domain=-3.15:-2.52,samples=20]
  ({etazx(\EtaRho,\t)},{etazy(\EtaRho,\t)});
\draw[fill=white,line width=.7pt]
  (1260,440) circle[radius=.057cm];
\end{scope}

\node[text=etaBlockBlue] at (1190,258) {$\delta_4$};
\node[text=etaBlockOrange] at (1220,619) {$\delta_4/2$};
\node[text=etaBlockBlue] at (1005,419) {$A_4$};
\draw[etaBlockBlue,thin] (1038,426)--(1126,445);
\node[text=etaBlockOrange] at (1120,802) {$C_4$};
\draw[etaBlockOrange,thin] (1141,783)--(1186,706);
\node[right] at (1275,440) {$1$};
\node at (1370,170) {$\mathbb T$};
\draw[etaBlockGray,thin] (1349,180)--(1255,215);

\draw[etaMain] (120,740)--(175,740);
\node[anchor=west] at (194,740) {main arc $A_4$};
\draw[etaConnector] (120,796)--(175,796);
\node[anchor=west] at (194,796) {circular connector $C_4$};
\draw[etaRadial] (120,852)--(175,852);
\node[anchor=west] at (194,852) {radial joins};
\node at (750,934)
  {$\operatorname{length}(A_4)\le2\varepsilon_4,
    \qquad\operatorname{length}(C_4)\le\varepsilon_4$};
\node at (750,991)
  {$\operatorname{length}(A_4)+\operatorname{length}(C_4)
    \le3\varepsilon_4$};
\end{tikzpicture}%
}
\caption{The two circular pieces at stage $n=4$.
The same geometry applies at every stage. The distances are measured
to $\mathbb T$.}
\label{fig:eta-circular-block}
\endgroup
\end{figure}

For $0<h<\varepsilon_1$, choose $N$ such that
\[
 \varepsilon_{N+1}<h\le\varepsilon_N.
\]
For $n<N$,
\[
 \frac{\delta_n}{2}
 >5\varepsilon_{n+1}
 \ge5\varepsilon_N
 \ge5h,
\]
so all circular pieces in these earlier groups miss
$B(\zeta,h)$.
For the group at $N$, the elementary
estimate
\[
 \HH(S\cap B(\zeta,h))\le2\pi h
\]
is valid 
for both circular arcs $S = A_N$ and $S=C_N$.  This shows that $A_N$ and $C_N$ contribute at most $4\pi h$. 
The remaining groups contribute at most
\[
 3\sum_{n>N}\varepsilon_n
 \lesssim\varepsilon_{N+1}
 \le h,
\]
up to an absolute multiplicative constant.
Together with the radial contribution, this gives
\[
 \HH(\eta\cap B(\zeta,h))\lesssim h.
\]
For $h\ge\varepsilon_1$, the same estimate follows from
\[
 \length(\eta)
 \le\delta_1+3\sum_{n=1}^{\infty}\varepsilon_n
 \lesssim\varepsilon_1.
\]
Consequently, we get $\sigma_\eta\in\CM(\mathbb D).$

\noindent\textbf{Step 4: Contradiction to \textup{$C_3$}.}
Set $\gamma=F(\eta)$. Since $F$ is conformal, this curve is simple and locally
rectifiable. Since it contains every $F(A_n)$,
\[
 \|\sigma_\gamma\|_{\CM(\Omega)}
 \ge
 \|\sigma_{F(A_n)}\|_{\CM(\Omega)}
 >n
 \qquad\text{for every }n.
\]
Hence $\sigma_\gamma\notin\CM(\Omega)$.

On the other hand, for every admissible map
$\phi:\Omega\to\Omega'$, the composition
$G=\phi\circ F$ maps $\mathbb D$ conformally onto a bounded
domain with rectifiable boundary. Bishop's
criterion~\eqref{eq:bishop}, applied to $\sigma_\eta$, gives
\[
 \length(\phi(\gamma))
 =
 \length(G(\eta))
 =
 \int_\eta|G'(z)|\,ds(z)
 <\infty.
\]
Thus, $\gamma$ does not satisfy \textup{$C_3$}. This  completes the proof.
\end{proof}

\section{Vanishing Carleson measures: Theorem~\ref{thm:vanishing-positive-measures}
and Proposition~\ref{prop:non}}
\label{sec:vanishing-measures}

We first show two preliminary facts on vanishing Carleson measures, which will be used in the later proof. 

Let \(U\) be a bounded domain and \(\mu\in\mathrm{CM}_0(U)\). We claim that \(\mu\) can be approximated in Carleson norm by its restrictions to compact subsets of \(U\). Choose an increasing sequence of compact sets \(K_j\subset U\) such that \(\bigcup_jK_j=U\).
Given \(\varepsilon>0\), the vanishing condition on $\mu$ provides \(r_0>0\) such that
\[
\frac{\mu(B(x,r)\cap U)}{r}<\varepsilon,
\qquad x\in\partial U,\quad 0<r<r_0.
\]
Since \(\mu\!\restriction_{U\setminus K_j}\le\mu\), the same bound holds for \(\mu\!\restriction_{U\setminus K_j}\) uniformly in \(j\). For \(r\ge r_0\), we instead have
\[
\frac{\mu(B(x,r)\cap(U\setminus K_j))}{r}
\le \frac{\mu(U\setminus K_j)}{r_0}.
\]
A Carleson measure on a bounded domain has finite total mass. Hence,  the approximation \(K_j\uparrow U\) implies that the right-hand side tends to zero as \(j\to\infty\). Combining these two estimates gives
\[
\bigl\|\mu-\mu\!\restriction_{K_j}\bigr\|_{\mathrm{CM}(U)}
\longrightarrow0.
\]
Each restriction \(\mu\!\restriction_{K_j}\) has compact support in \(U\) and therefore belongs to \(\mathrm{CM}_0(U)\).

For a positive Borel measure $\nu$ on $\mathbb D$, put
\[
 J_\nu(a)
 :=
 \int_{\mathbb D}
 \frac{1-|a|^2}{|1-\overline a z|^2}\,d\nu(z),
 \qquad a\in\mathbb D.
\]
The classical disk kernel criterion states that
\begin{equation}\label{eq:vanishing-kernel}
 \nu\in\CM_0(\mathbb D)
 \quad\Longleftrightarrow\quad
 \lim_{|a|\to1}J_\nu(a)=0
\end{equation}
where the limit is uniform in the argument of $a$; 
see, for example, \cite[Theorem~8.2.5]{Zhu1990}.

\begin{proof}[Proof of Theorem~\ref{thm:vanishing-positive-measures}] 
We establish the cycle of implications
\[
\textup{$V_1$}
\Longrightarrow
\textup{$V_2$}
\Longrightarrow
\textup{$(\ast)$}
\Longrightarrow
\textup{$V_1$},
\]
where \textup{$(\ast)$} is the auxiliary assertion that
$(F^{-1})^*$ maps every  measure in
$\CM_0(\mathbb D)$ to a measure in $\CM_0(\Omega)$. 
The implication from Ahlfors \(1\)-regularity of \(\Gamma\) to \(\textup{$(\ast)$}\) is already known; see \cite{WeiZinsmeister}. To the best of our knowledge, the converse has not previously been recorded.

\noindent\textbf{\textup{$V_1$} implies \textup{$V_2$}.}
Under Ahlfors $1$-regularity, applying \eqref{eq:transfer} gives 
$\mu\in\CM_0(\Omega)$ if and only if $F^*\mu\in\CM_0(\mathbb D)$. Indeed, 
 since the map $F$ preserves compact
support, by compact approximation explained above and
closedness of vanishing Carleson measures, we conclude that 
\[
 \mu\in\CM_0(\Omega)
 \quad\Longleftrightarrow\quad
 F^*\mu\in\CM_0(\mathbb D).
\]
Writing $\nu=F^*\mu$, we have 
\[
 \int_\Omega|\phi_a'(w)|\,d\mu(w)
 =
 \int_{\mathbb D}|m_a'(z)|\,d\nu(z)
 =
 J_\nu(a).
\]
Applying \eqref{eq:vanishing-kernel} proves \textup{$V_2$}.

\noindent\textbf{\textup{$V_2$} implies \textup{$(\ast)$}.}
For a measure $\nu\in\CM_0(\mathbb D)$, put
$\mu=(F^{-1})^*\nu$. Then $F^*\mu=\nu$, so
\[
 \int_\Omega|\phi_a'(w)|\,d\mu(w)
 =
 J_\nu(a)\longrightarrow0.
\]
Condition \textup{$V_2$} yields
$\mu\in\CM_0(\Omega)$.

\noindent\textbf{\textup{$(\ast)$} implies \textup{$V_1$}.}
First, \textup{$(\ast)$} yields the uniform bound
\begin{equation}\label{eq:vanishing-transfer-bound}
 \|(F^{-1})^*\nu\|_{\CM(\Omega)}
 \le C\|\nu\|_{\CM(\mathbb D)},
 \qquad \nu\in\CM_0(\mathbb D).
\end{equation}
Otherwise, by rescaling positive measures, we could choose
$\nu_n\in\CM_0(\mathbb D)$ with
\[
 \|\nu_n\|_{\CM(\mathbb D)}\le2^{-n},
 \qquad
 \|(F^{-1})^*\nu_n\|_{\CM(\Omega)}\ge n.
\]
The positive sum $\nu=\sum_{n=1}^\infty\nu_n$ converges
in Carleson norm and belongs to $\CM_0(\mathbb D)$.
Positivity gives
\[
 \|(F^{-1})^*\nu\|_{\CM(\Omega)}
 \ge
 \|(F^{-1})^*\nu_n\|_{\CM(\Omega)}
 \ge n
 \qquad\text{for every }n,
\]
contradicting \textup{$(\ast)$}.

Now take any  $\nu\in\CM(\mathbb D)$.
Its compact restrictions
\[
 \nu_k=\nu\!\restriction_{\{|z|\le1-1/k\}},
 \qquad k\ge2,
\]
belong to $\CM_0(\mathbb D)$ and satisfy
$\|\nu_k\|_{\CM(\mathbb D)}\le\|\nu\|_{\CM(\mathbb D)}$.
Set
\[
 \mu_k=(F^{-1})^*\nu_k,
 \qquad
 \mu=(F^{-1})^*\nu.
\]
By \eqref{eq:pushforward} and monotone convergence,
$\mu_k(E)\uparrow\mu(E)$ for every Borel set $E\subset\Omega$.
The uniform bound \eqref{eq:vanishing-transfer-bound}
therefore gives
\[
 \mu(B(x,r)\cap\Omega)
 =
 \lim_{k\to\infty}\mu_k(B(x,r)\cap\Omega)
 \le C\|\nu\|_{\CM(\mathbb D)}\,r
\]
for all $x\in\Gamma$ and $r>0$.
Thus, $\mu\in\CM(\Omega)$. The  characterization
\eqref{eq:push-characterization} now implies
\textup{$V_1$}.
\end{proof}

\medskip

\begin{proof}[Proof of Proposition~\ref{prop:non}]
We first show that, for every bounded Jordan domain $\Omega$
and every simple locally rectifiable curve
$\gamma\subset\Omega$,
\begin{equation}\label{eq:vanishing-curve}
 \sigma_\gamma\in\CM_0(\Omega)
 \quad\Longleftrightarrow\quad
 \length(\gamma)<\infty
 \ \text{and}\
 \dist(\gamma,\Gamma)>0.
\end{equation}
Suppose $\sigma_\gamma\in\CM_0(\Omega)$.
Its total mass is finite, so $\gamma$ has finite length.
If $\dist(\gamma,\Gamma)=0$, boundedness of $\Omega$
would give an accumulation point $x\in\Gamma$
of $\gamma$. Since $\gamma$ is connected, the $1$-Lipschitz
map $z\mapsto|z-x|$ sends $\gamma$ onto an interval with
infimum zero and containing a positive number.
It therefore sends $\gamma\cap B(x,r)$ onto a set
containing $(0,r)$ for every sufficiently small $r>0$.
Consequently,
\[
 \sigma_\gamma(B(x,r))
 =
 \HH(\gamma\cap B(x,r))
 \ge r,
\]
contrary to the vanishing condition.

Conversely, suppose $\gamma$ has finite length and
$d:=\dist(\gamma,\Gamma)>0$.
Boundary disks of radius less than $d$ miss $\gamma$,
whereas for $r\ge d$,
\[
 \frac{\sigma_\gamma(B(x,r))}{r}
 \le\frac{\length(\gamma)}{d}.
\]
Thus, $\sigma_\gamma\in\CM_0(\Omega)$, proving
\eqref{eq:vanishing-curve}.

The conjunction of finite length and positive boundary
distance appeared in \eqref{eq:vanishing-curve} is conformally invariant between bounded domains.
Indeed, positive boundary distance is equivalent to relative
compactness, which is preserved by a conformal bijection.
On the corresponding compact subsets, the map and its
inverse have bounded derivatives, so finiteness of length
is preserved as well. With $\eta=F^{-1}(\gamma)$,
\eqref{eq:vanishing-curve} therefore gives
\[
 \sigma_\gamma\in\CM_0(\Omega)
 \quad\Longleftrightarrow\quad
 \sigma_\eta\in\CM_0(\mathbb D).
\]
Finally, the arclength change-of-variables formula yields
\[
 \int_\gamma|\phi_a'(w)|\,ds(w)
 =
 \int_\eta|m_a'(z)|\,ds(z)
 =
 J_{\sigma_\eta}(a).
\]
The  criterion \eqref{eq:vanishing-kernel} completes
the proof.
\end{proof}

\paragraph*{Acknowledgments.}
Research is supported  by the National Natural Science Foundation
of China (grant nos.~12271218 and 12571083).
\vspace{-0.3cm}
\paragraph*{AI disclosure.}
The authors used GPT-6.1 Sol Ultra to assist with the preparation of the figure
and to improve the language and readability of the manuscript.

\end{document}